\documentclass[11pt]{article}
\usepackage{epsf}
\usepackage{color}
\usepackage{epsfig}

\usepackage{graphicx}
\usepackage{xcolor}
\usepackage{amsmath}
\usepackage{latexsym}
\usepackage{amssymb}
\usepackage{amsfonts}

\usepackage{a4}
\usepackage{longtable}

\renewcommand{\H}{H}
\newcommand{\Emin}{E_{\rm min}}
\newcommand{\Emax}{E_{\rm max}}

\newcommand{\e}{\ensuremath{\mathrm{e}}}

\renewcommand{\P}{P}
\newcommand{\R}{\mathbb{R}}
\newcommand{\C}{\mathbb{C}}
\newcommand{\T}{t}

\newcommand{\structurecomment}[1]{}

\begin{document}

\title{An efficient algorithm based on splitting for the time integration of the Schr\"odinger equation}

\author{Sergio Blanes$^{1}$\thanks{Email: \texttt{serblaza@imm.upv.es}}
   \and
  Fernando Casas$^{2}$\thanks{Email: \texttt{Fernando.Casas@uji.es}}
   \and
Ander Murua$^{3}$\thanks{Email: \texttt{Ander.Murua@ehu.es}}
   }

\maketitle

\begin{abstract}

We present a practical algorithm based on symplectic splitting methods to integrate
numerically in time the Schr\"odinger equation. When
discretized in space, the Schr\"odinger equation can be recast as a
classical Hamiltonian system 
corresponding to a generalized high-dimensional separable harmonic
oscillator. The particular structure of this  system combined with previously obtained stability
and error analyses
allows us to construct a set of highly efficient symplectic
integrators  with sharp error bounds and optimized for different
tolerances and time integration intervals. They can be considered, in this
setting, as polynomial approximations to the matrix exponential in a similar way as methods based on
Chebyshev and Taylor polynomials.
 The theoretical analysis, supported by
numerical experiments, indicates that the new methods are more efficient than
schemes based on Chebyshev polynomials for all
tolerances and time intervals. The algorithm we present incorporates
the new splitting methods and automatically selects the most efficient scheme
given a tolerance, a time integration interval and an estimate on the spectral radius of the Hamiltonian.

\vspace*{0.5cm}

\begin{description}
 \item $^1$Instituto de Matem\'atica Multidisciplinar,
  Universitat Polit\`ecnica de Val\`encia, E-46022  Valencia, Spain.
 \item $^2$Institut de Matem\`atiques i Aplicacions de Castell\'o and
   Departament de Matem\`atiques, Universitat Jaume I,
  E-12071 Castell\'on, Spain.
 \item $^3$Konputazio Zientziak eta A.A. saila, Informatika
Fakultatea, UPV/EHU, Donostia/San Sebasti\'an, Spain.
\end{description}

\end{abstract}

\section{Introduction}

When investigating the dynamical behavior of quantum systems of low to moderate dimension,
very often it is necessary to
solve numerically the time dependent Schr\"odinger equation ($\hbar = 1$)
\begin{equation}   \label{Schr0}
  i \hbar \frac{\partial}{\partial t} \psi (x,t)  = \hat{H} \psi (x,t), \qquad \psi(x,0)=\psi_0(x).
\end{equation}
Here $\hat{H}$ is the Hamiltonian operator, $\psi:\mathbb{R}^d\times\mathbb{R}
\longrightarrow \mathbb{C}$ is the wave function representing the
state of the system and $\psi_0(x)$ is the initial state. For simplicity, in the sequel we consider $\hat{H} = \hat{T} + \hat{V}$, with the kinetic energy operator
$\hat{T} = -\Delta/(2 \mu)$ for a reduced mass $\mu>0$ and a potential $\hat{V}$, although the procedure presented in
this paper is also valid for more general Hamiltonian operators.

The solution of (\ref{Schr0}) can be expressed as
\begin{equation}  \label{evol1}
  \psi(x,t) = \hat{U}(t) \psi_0(x),
\end{equation}
the (unitary) evolution operator $\hat{U}$ being formally given by
  $ \hat{U}(t) = \e^{-i t \hat{H}}$. In practice, however, it is not possible to get a closed expression for $\hat{U}(t)$,
  and so numerical methods are applied to get reliable approximations. This process
involves typically two stages. In the first a discrete spatial representation of the  initial wave function $\psi_0(x)$ and the operator $\hat{H}$
on an appropriate grid are constructed. In the second, this finite representation is propagated in time
with a numerical integrator.

As for the space discretization process, several techniques can be used, depending on the particular
problem one aims to analyze: finite difference schemes,
spectral methods based on collocation with trigonometric polynomials, Galerkin method with a Hermite basis,
etc, both in one or more dimensions (see \cite{lubich08fqt} and references therein).
The space discretization process restricts the energy range of the approximation
and imposes an upper bound to the high frequency components
represented by the discrete solution.

In any event, once this process has been carried out,
one has the linear system of ordinary differential equations 
\begin{equation} \label{td.1}
  i \frac{d }{dt} u(t) = H u(t), \qquad u(0)=u_0 \in \mathbb{C}^N,
 \end{equation}
where $u(t)$ now represents a discretized version of the wave function $\psi(x,t)$ at the $N$ space grid points, with $N$ usually a large number.
 The goal is then to compute $u(t)$ at a given target time $\T$ from the known value of $u(0)=u_0$.
The $N \times N$ matrix $H$ (and in particular its
discrete spectrum) depends of course on the particular space discretization carried out.
We will hereafter assume that $H$ is a real symmetric matrix which implies that it can be diagonalized
with real eigenvalues.

The exact solution of eq. (\ref{td.1}) reads
 \begin{equation}   \label{td.2}
    u(t) = \e^{-i \, t \, H} \, u_0,
\end{equation}
but computing the matrix exponential $\e^{-i\,  \T\,  H}$
by diagonalizing $H$ (usually, a matrix of large dimension and large norm) is prohibitively expensive.
An effective alternative consists in computing approximations of $u(\T)$ of the form
\begin{equation}   \label{approx.1}
  u(\T) \approx \P_m(\T\, H) u_0,
\end{equation}
where $\P_m(y)$ is a polynomial in $y$ that approximates the exponential $\e^{-i\, y}$, since in that case only multiplications of the
matrix $H$ with vectors $u$ are necessary. These products can be efficiently evaluated  in
complex variables  (provided that  a Fourier spectral method is used to obtain the discretized version (\ref{td.1})  of (\ref{Schr0})) with the complex-to-complex 
Fast Fourier Transform (FFT) algorithm \cite{feit82sot,kosloff83afm,kosloff88tdq,leforestier91aco}.

There are different choices for such a polynomial $\P_m(y)$. For instance, one may consider truncated Taylor or Chebyshev series expansion of $\e^{-i\, y}$
for an appropriate real interval of $y$, or the Lanczos method, where the polynomial is determined by a Galerkin approximation on the Krylov space
spanned by $u_0, H u_0, \ldots, H^{m-1} u_0$ \cite{park86uqt}.

In this paper we consider yet another kind of polynomial approximation to $\e^{-i\, \T\, H}u_0$, namely one based on explicit symplectic splitting methods
\cite{gray96sit,gray94chs,blanes06sso,blanes08otl,blanes11eao}.
This approach can be applied under the same assumptions than the Chebyshev method, the main difference being the following. Whereas in the
Chebyshev (or Taylor) method the approximation (\ref{approx.1}) is constructed by evaluating products of the form $H u$, where $u \in \mathbb{C}^N$,
with symplectic splitting methods one writes $u = q + i p$, $q,p \in \mathbb{R}^N$. The algorithm then proceeds by successively computing \emph{real} matrix-vector
products $H q$ and $H p$ with different weights, so that the real and imaginary parts of $\e^{-i \, \T\,  H}u_0$ are approximated in a different way, with a
much reduced computational cost.

More specifically, if a spatial discretization based on Fourier spectral methods is considered, then the cost of computing $H u$, $u \in \mathbb{C}^N$,
 amounts essentially to
one complex-to-complex FFT and its inverse, whereas in the case of $H v$, $v \in \mathbb{R}^N$, one has to evaluate one real-to-complex FFT and its
inverse complex-to-real FFT, and this process requires half the computing time of the fully complex case. As a result, the proposed
algorithm based on splitting methods
turns out to be between $1.4$ and 2 times faster than the Chebyshev method for the same accuracy in all the examples 
we have analyzed.
Moreover, the procedure is easy to implement and the resulting approximations
preserve important qualitative properties of the
exact solution.

The algorithm we present here has embedded
several symplectic splitting schemes designed according to different optimization criteria with the purpose of covering most of the cases one
finds in practical applications (high accuracy over long time intervals,
low accuracy over short times, etc.).
The computation of the coefficients of the methods, which constitutes a non-trivial task by itself,
is largely based on the stability and error analysis of splitting methods carried out in \cite{blanes08otl,blanes11eao}.
Given a target value of time $\T$ and an error tolerance, the algorithm selects a specific symplectic splitting scheme leading to a numerical solution with
the prescribed accuracy and the minimum computational work, measured as the number of real matrix-vector products. By construction,
the algorithm developed here is aimed to be applied for the same problems and under the same assumptions as the Chebyshev method, with a remarkable
gain in efficiency for all the examples we have tested.

The plan of the paper is the following. Since our procedure may be considered as an alternative to the Chebyshev method,
in section~\ref{sec.2} we
summarize the main features of the schemes based on
this polynomial approximation of $\e^{-i t H} u_0$.  In section~\ref{sec.3} we
analyze the stability and the global error of symplectic splitting methods in this context, and the actual algorithm is presented, whereas the comparison with
Chebyshev (and Taylor as a reference) is carried out in section~\ref{sec.5} on a pair of selected numerical examples. 

\section{Polynomial approximations}
 \label{sec.2}

\subsection{General considerations}
\label{subsec.2.1}

Given a $m$th degree polynomial $P_m(y)$ approximating $\e^{-i\, y}$,
the solution $u(t)=\e^{-i \, t \, H}u_0$ of (\ref{td.1}) at a prescribed target time $\T$ can be approximated as
\begin{equation}\label{eq:w_tau}
u(\T) \approx  u_1= P_m(\T \, \H)\, u_0,
\end{equation}
with the corresponding error (in Euclidean norm)  bounded as
\begin{equation*}
    \|u_1 - \e^{-i\, \T\, \H}u_0\| \leq  \max_{j=0,1,\ldots, N-1} |\P_m(\T\, E_j) -\e^{-i\, \T\, E_j}| \,\|u_0\|
\end{equation*}
in terms of the (real) eigenvalues $E_0,\ldots,E_{N-1}$ of $\H$. Assuming that
the spectrum $\sigma(H)=\{E_0,\ldots,E_{N-1}\}$ is contained in an interval of the form $[\Emin,\Emax]$, then
\begin{equation*}
    \|u_1 - \e^{-i\, \T\, \H}u_0 \| \leq  \sup_{\T \, \Emin \leq y \leq \T\, \Emax} |\P_m(y) - \e^{-i\, y}|\, \|u_0\|.
\end{equation*}

There are several possibilities to estimate $\Emax$ and
$\Emin$ for different classes of matrices  (see e.g. \cite{huang05ase,mazzi11drf,yang11asm,zhu11ets}). If $H$ can be decomposed as the sum $H = T + V$ of two symmetric matrices with known lower and upper bounds for their eigenvalues, $\Emin$ (resp. $\Emax$) can be simply obtained as the sum of the lower (resp. upper) bounds of the eigenvalues of $T$ and $V$. This happens, in particular, when the Hamiltonian operator $\hat{H}= -\Delta/(2 \mu) + \hat{V}$  is discretized by spectral Fourier collocation with $N$ Fourier modes, in which case
\begin{equation} \label{eq:bounds}
   \Emin = \min_{x} V(x), \qquad \Emax = \frac{1}{2 \mu} \frac{N^2}{4} + \max_x V(x).
\end{equation}
In any event, once $\Emin$ and $\Emax$ have been determined, we introduce
\begin{equation} \label{shifting}
\alpha = \frac{\Emax + \Emin}{2},
\qquad \beta = \frac{\Emax - \Emin}{2}, \qquad
   \overline{H} = H - \alpha I,
\end{equation}
so that the spectrum of the shifted operator $\overline{H}$ is contained in an interval centered at the origin,
$\sigma(\overline{H}) = \{E_0-\alpha,\ldots,E_{N-1}-\alpha\} \subset [-\beta,\beta]$. We thus have
 \begin{equation}   \label{td.2Cheb}
    \e^{-i \, \T \,  \H} \, u_0 = \e^{-i\, \T\,  \alpha} \, \e^{-i\, \T\, \overline{H}} \, u_0.
\end{equation}
Hence, we will hereafter assume without loss of generality that our problem consists in approximating $\e^{-i\,\T\,\H}u_0$
for a real symmetric matrix $\H$ with $\sigma(\H) \subset [-\beta,\beta]$. In that case,
\begin{equation}
\label{eq:Err1}
  \|u_1 - \e^{-i\, \T\, \H}u_0 \| \leq  \epsilon_m(\beta\, \T)\, \|u_0\|,
\end{equation}
where
\begin{equation}
  \label{eq:Err2}
\epsilon_m(\theta) \equiv  \sup_{-\theta\leq y \leq \theta}|\e^{-i \, y} -
P_m (y) |.
\end{equation}

\subsection{Taylor polynomial approximation}
 \label{sec.2.1}

The $m$th degree Taylor polynomial $P_m^T(y)$ corresponding to $\e^{-i\, y}$ is of course
\begin{equation}\label{eq:Taylor}
 P_m^T(y) \equiv \sum_{k=0}^m\frac{(-i)^k}{k!} y^k,
\end{equation}
and Horner's algorithm provides an efficient way to compute $u_1 =  P_m^T(\T\, \H) u_0$, namely
\begin{equation}\label{alg:Horner}
\begin{array}{l}
  {y}_0  = {u}_0 \\
  {\bf do} \ \ k=1,m \\
     \quad {y}_k  = {u}_0 - i \displaystyle \frac{\T}{m+1-k}  \H {y}_{k-1} \\
  {\bf enddo} \\
u_1 = y_m.
\end{array}
\end{equation}
The process requires storing three complex vectors (or
equivalently, 6 real vectors).

An error estimate of the form
(\ref{eq:Err1}) can be obtained with $\epsilon_m(\theta)$ in (\ref{eq:Err2}) replaced by its upper bound
\begin{equation}\label{eq:ErrorTaylor}
 \epsilon_m^T(\theta)  \equiv \frac{\theta^{m+1}}{(m+1)!}.
\end{equation}
Since $m! \sim \sqrt{2\pi m} \ (m \e)^{m}$ for large values of $m$ \cite{olver10nho}, we can write
\[
    \epsilon_m^T(\theta)\sim \frac{1}{\e\, \sqrt{2\pi m}}
\left(\frac{\theta\, \e}{m}\right)^{m+1}.
\]
In consequence,  we cannot expect to have a reasonably accurate
approximation $P_m^T(\T\, \H) u_0$ of $\e^{-i \, \T\, \H}u_0$ unless
\[
    m > \e \, \theta = \e\, \beta\, t.
\]
In other words, increasing the value of the target time $\T$ where the solution is to be found and/or refining the spatial discretization (so that
$\beta$ gets larger) requires evaluating a higher degree Taylor polynomial.




%
\subsection{Chebyshev polynomial approximation}
 \label{sec.2.2}

The Chebyshev polynomial expansion scheme, proposed for the first
time in the context of the Schr\"odinger equation  in
\cite{tal-ezer84aaa}, constitutes a standard tool to compute
(\ref{td.2}). A detailed analysis of the
procedure, including error estimates for the problem at hand,
can be found in \cite{lubich08fqt}. For completeness, we review here some of
its main features.

The $m$th degree truncation of the Chebyshev series expansion of
$\e^{-i\, y}$ in the interval $y \in [-\theta,\theta]$ is given by
\begin{equation}\label{eq:Chebyshev}
 P_{m,\theta}^C(y) \equiv    J_0(\theta) + 2 \sum_{k=1}^{m} (-i)^k J_k(\theta) \, T_k(y/\theta) ,
\end{equation}
where for each $k$, $J_k(t)$ is the Bessel
function of the first kind \cite{olver10nho} and $T_k(x)$ is the $k$th Chebyshev
polynomial generated from the recursion
\begin{equation}\label{ChebPolyn}
  T_{k+1}(x)=2xT_k(x)-T_{k-1}(x), \qquad k\geq 1
\end{equation}
and $T_0(x)=1, \ T_1(x)=x$. According with the analysis in \cite{lubich08fqt},
$\e^{-i\, \T\, \H}u_0$ can be approximated by $P_{m,\beta\T}^C(\T\, \H) u_0$
 with an error estimate of the form (\ref{eq:Err1}), where $ \epsilon_m(\theta)$ in (\ref{eq:Err2}) is replaced by its upper bound
\begin{equation}\label{eq:ErrorChebyshev}
 \epsilon_m^C(\theta)  \equiv
4 \left( \e^{1-\theta^2/(2m+2)^2} \frac{\theta}{2m+2}
  \right)^{m+1}.
\end{equation}
%

\begin{figure}
\begin{center}
  \includegraphics[width = 1\textwidth]{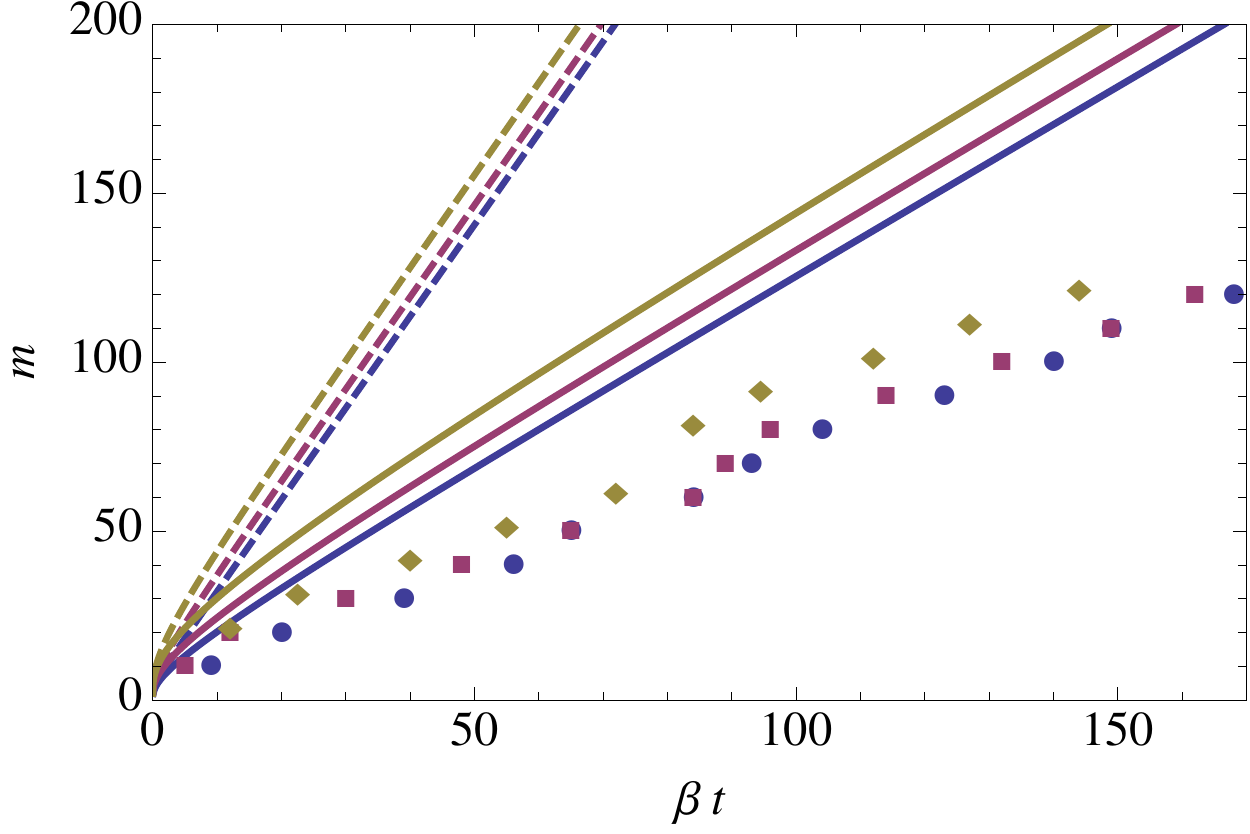}
\caption{\label{fig:mTayChev} Comparison of the required minimum polynomial degree $m$ as function of $\theta=\beta\, t$ for Taylor (dashed line) and
Chebyshev (continuous line) for different values of error tolerance: $\texttt{tol}=10^{-4}, 2 \times 10^{-7}, 10^{-11}$. 
Diamonds, squares and circles stand for the
computational cost (equivalent to a polynomial approximation of degree $m$) 
for error tolerances below $10^{-4}$, $2 \times 10^{-7}$ and $10^{-11}$, respectively, obtained with
 symplectic splitting schemes in Table~\ref{tab:relevant_parameters}.}
\end{center}
\end{figure}

 In Figure~\ref{fig:mTayChev} we depict  the minimum degree $m$ as a function of $\theta=\beta\, t$ of Chebyshev approximations
 for prescribed tolerances $\texttt{tol}=10^{-4}, 2 \times 10^{-7}, 10^{-11}$, so that $\epsilon_m^C(\beta\, t)\leq \texttt{tol}$ (continuous lines)
 in comparison with the corresponding degree $m$ for Taylor approximations
 (dashed lines) such that $\epsilon_m^T(\beta\, t) \leq \texttt{tol}$. Notice that Chebyshev always gives a similar accuracy with a lower degree polynomial
 (hence, with less computational cost), with a gain in efficiency of up to a factor of two for sufficiently large values of $\theta = \beta t$.

Once the degree of the polynomial $m$ has been chosen, given a certain error tolerance, target time $t$, and bound $\beta$ of $\sigma(\H)$,
one has to compute $\P_{m,\beta\T}^C(\T \, \H)\, u_0$ in an as efficient as possible way. This can be done with
the Clenshaw recursive algorithm as follows: first evaluate
 the coefficients $c_k=(-1)^k J_k(\beta t)$ for $k=0,1,\ldots,m$ and then compute recursively
\begin{equation}
 \label{alg:Clenshaw}
\begin{array}{l}
  {d}_{m+2}  = {0}, \qquad   {d}_{m+1}  = {0} \\
  {\bf do} \ \ k=m,m-1,\ldots,1,0 \\
     \quad {d}_k  = c_k\, {u}_0 + \frac{2}{\beta} \H {d}_{k+1} - {d}_{k+2} \\
  {\bf enddo} \\
u_1  = {d}_0 - {d}_2,
\end{array}
\end{equation}
which produces $u_1 \equiv P^C_{m,\beta t}( t \H) \, u_0 \approx \e^{-i\, t\, H} u_0$ as output.
Clenshaw algorithm keeps only four complex
vectors in memory\footnote{If the vectors are written in their real and
imaginary part, and the algorithm is carried out in real
variables, then the algorithm needs to store only seven real
vectors instead of eight.}, but the whole procedure has to be carried out for each value of $m$.
Since the coefficients $c_k$ are relatively small
as $k$ grows, the Clenshaw algorithm is stable and so it is possible to work with polynomials of very high degree
(even in the thousands) provided the Bessel functions are accurately computed.

\section{Symplectic splitting methods}
 \label{sec.3}

\subsection{General considerations}

An alternative to Chebyshev polynomial approximations of $\e^{-i\, t\,
  \H}u_0$ first considered in \cite{gray96sit,gray94chs} consists in
applying specially designed splitting methods to numerically
integrate the system (\ref{td.1}) recast in a more suitable form.

By considering $q = \mbox{Re}(u) \in \R^N$ and $p = \mbox{Im}(u) \in \R^N$,
 equation (\ref{td.1}) is equivalent to
\begin{equation} \label{td.1b}
  \frac{d }{dt} z = (A+B) z,\quad z(0)=z_0,
 \end{equation}
where
\begin{equation}    \label{m.a.b}
  z = \left(  \begin{array}{c}
               q \\
               p
              \end{array} \right),
              \qquad
  A = \left(  \begin{array}{ccc}
               0   & \, &  H  \\
               0    & \, &  0  \end{array} \right),
              \qquad
  B = \left(  \begin{array}{ccc}
               0   & \, & 0  \\
              -H    & \, & 0  \end{array} \right).
\end{equation}
The solution $z(t) = \e^{t \, (A+B)}z_0$
of (\ref{td.1b}) can be written 
in terms of the orthogonal and symplectic matrix
\begin{equation}  \label{O(y)}
   O(y) = \left(
 \begin{array}{rcr}
   \cos(y)   &  & \sin(y)  \\
  -\sin(y)  &   & \cos(y)  \end{array}
  \right)
\end{equation}
as $z(t) = O(t \, H) z_0$.
To introduce general symplectic splitting methods in this setting, let us first show how the well known Strang splitting
can
be used to approximate $\e^{-i\, t\,  \H}u_0$.  Let $m$ be a sufficiently large positive integer, so
that for $\tau=t/m$, we consider the approximation
\[
  \e^{\tau (A+B)} \approx \e^{\frac{\tau}{2}A} \, \e^{\tau B} \,  \e^{\frac{\tau}{2}A}.
\]
It is then clear that
\[
 \e^{t \, (A+B)} = \left( \e^{\tau\, (A+B)}\right)^m \approx
\left( \e^{\frac{\tau}{2}A} \, \e^{\tau B} \,
\e^{\frac{\tau}{2}A}\right)^m =
\e^{\frac{\tau}{2}A} \left( \e^{\tau B} \,
\e^{\tau  A}\right)^{m-1} \, \e^{\tau B} \,\e^{\frac{\tau}{2}A},
\]
or equivalently,
\begin{equation} \label{compos}
 O(t\, \H) = \e^{t\, (A+B)}  \approx   K(t\,  \H) = \e^{t\,  a_{m+1}\, A} \,\e^{t\, b_{m}\, B} \,\e^{t\, a_m\, A} \ \cdots \
  \e^{t\, b_1\, B} \, \e^{t\, a_1\, A },
\end{equation}
with
\begin{equation}
  \label{eq:abStrang}
  (a_1,b_1,a_2,\ldots,a_{m},b_{m},a_{m+1}) = \left(\frac{1}{2m},\frac{1}{m},\frac{1}{m},\ldots,\frac{1}{m},\frac{1}{m},\frac{1}{2m}\right).
\end{equation}
Due to the nilpotent structure of the matrices $A$ and $B$ in
(\ref{m.a.b}), the exponentials in the definition (\ref{compos}) of $K(\T\, \H)$ take a particularly simple form, namely
\begin{equation}    \label{exp.a.b}
  \e^{t\,  a_j\,  A} = \left(
 \begin{array}{ccc}
    I   & \,  & a_j\,  t\, \H  \\
    0  & \,  &  I   \end{array} \right), \qquad \qquad
   \e^{t\,  b_j\,  B} = \left(
 \begin{array}{ccc}
   I   & \, & 0 \\
  -b_j\,  t\,  \H  &  \, & I  \end{array} \right).
\end{equation}
This analysis shows that the approximation $K(\T\, \H)z_0 \approx
\e^{\T\, (A+B)}z_0$
can be computed with the following procedure, similar in nature and equivalent in computing time to the
Horner (\ref{alg:Horner}) and Clenshaw  (\ref{alg:Clenshaw}) algorithms: Given $u_0\in \C^N$,
\begin{equation}
  \label{alg:splitting}
  \begin{array}{l}
 q := \mbox{Re}(u_0),\\
 p := \mbox{Im}(u_0),\\
  {\bf do} \ \ k=1,m \\
     \quad q  := q + a_k\,  \T \, \H\,p \\
     \quad p  := p - b_k\,  \T\,  \H\,q \\
  {\bf enddo} \\
     \quad q  := q + a_{m+1} \, \T\,  \H\,p \\
     \quad u_1 := q + i p,\\
\end{array}
\end{equation}
producing $u_1 \approx \e^{-i\, \T\, \H}u_0$ as output.
Notice that it only requires storing three real vectors of dimension
$N$ (namely $q$, $p$, and $w=Hp$ or $w=Hq$) instead of seven real vectors for the Clenshaw algorithm and six real vectors for the Horner algorithm.
It is worth remarking that, since $\e^{A}$ and $\e^B$ are symplectic matrices, $K(t H)$ is also symplectic. Unitarity is
no longer preserved by this scheme, but neither the average error in energy nor the norm of the solution increases with time, since it is conjugate to
a unitary method \cite{blanes08otl}.

In practice, and in the same way as other polynomial approximations, it is convenient to apply Algorithm (\ref{alg:splitting})
with the original $H$ replaced by the shifted version $\overline{H}$ considered in (\ref{shifting}) (and then make use of the equality (\ref{td.2Cheb})), so that the spectrum of $\overline{H}$ is contained in an interval of the form $[-\beta,\beta]$ with $\beta$ as sharp as possible.
Therefore, in what follows we always
assume that $\sigma(H) \subset [-\beta,\beta]$.

Although Algoritm (\ref{alg:splitting}) with coefficients (\ref{eq:abStrang}) can be used in principle to approximate $\e^{-i\, \T\,\H}u_0$,
we next show that, for given values of $m$ and
$\theta = \beta t$, much better approximations can be obtained if other sequences of coefficients
$(a_1,b_1,a_2,\ldots,a_{m},b_{m},a_{m+1})$
are chosen instead. To see how this can be done, an error estimate of the corresponding approximation (\ref{compos}) is necessary first.

\subsection{Error analysis}
\label{sec.3.2}

For a given finite sequence of real numbers
\begin{equation}
  \label{eq:sequence}
  (a_1,b_1,a_2,\ldots,a_{m},b_{m},a_{m+1}),
\end{equation}
Algorithm (\ref{alg:splitting}) produces an approximation of the form
\[
\left(
  \begin{matrix}
    q_1 \\ p_1
  \end{matrix}
\right)
=K(\T\, \H)  \left(
  \begin{matrix}
    q_0 \\ p_0
  \end{matrix}
\right) \approx \e^{\T\, (A+B)} \left(
  \begin{matrix}
    q_0 \\ p_0
  \end{matrix}
\right)
\]
(or equivalently, $q_1 + i \, p_1 \approx \e^{-i\,\T\, \H} (q_0+i\, p_0)$)
with
\begin{equation}\label{K(x)}
  K(\T\, \H) =  \left(
 \begin{array}{cc}
  K_{11}(\T\, \H) &  K_{12}(\T\, \H) \\
  K_{21}(\T\, \H) &  K_{22}(\T\, \H)
 \end{array} \right).
\end{equation}
Here $K_{11}(y)$, $K_{22}(y)$ are even polynomials of degree $2m$, $K_{12}(y)$ and $K_{21}(y)$ are odd polynomials of degree
$2m-1$ and $2m+1$ respectively, and
$\det K(y) = K_{11}(y) K_{22}(y) - K_{12}(y) K_{21}(y) \equiv 1$.
It is important to remark that for a given positive integer $m$, compared to Horner's  (\ref{alg:Horner}) and Clenshaw's (\ref{alg:Clenshaw}) algorithms,
the degree of the polynomials
involved in an $m$-stage splitting method (\ref{eq:sequence}) is \emph{twice} the degree of the corresponding Taylor and Chebyshev
polynomials, with the same computational cost.

\subsubsection{Error estimates for a single application of a splitting method}
\label{sss:err}

We next focus on obtaining upper bounds for the error
\begin{eqnarray*}
  \|(q_1 + i \, p_1) -\e^{-i\,\T\, \H} (q_0+i\, p_0) \| &=&
\left\|
K(\T\, \H) \left(
  \begin{matrix}
    q_0 \\ p_0
  \end{matrix}
\right)
- O(\T\, \H) \left(
  \begin{matrix}
    q_0 \\ p_0
  \end{matrix}
\right)
\right\| \\
&\leq &
\left\|
K(\T\, \H)
- O(\T\, \H)
\right\|
\, \|q_0+ i\, p_0 \|
\end{eqnarray*}
in Euclidean norm. Since $H$ is assumed to be a real symmetric matrix, it can be diagonalized as
\begin{equation*}
  H = P^T \,
\left(
  \begin{matrix}
    E_0 & 0 &  \cdots & 0 \\
    0   & E_1 &  \cdots & 0 \\
    0   &  0  & \ddots &  0 \\
    0   & \cdots & 0 & E_{N-1}
  \end{matrix}
\right) \, P,
\end{equation*}
where $P$ is an orthogonal $N\times N$ matrix.  We thus have
\begin{equation*}
K(\T\, \H)
- O(\T\, \H) =
P^T \,
\mathcal{E}
\, P,
\end{equation*}
where $ \mathcal{E}$ is the block-diagonal matrix (with $2\times 2$ matrices at the diagonal)
\begin{equation*}
\left(
  \begin{matrix}
    K(\T\, E_0)-O(\T\, E_0) & 0 &  \cdots & 0 \\
    0   &  K(\T\, E_1)-O(\T\, E_1) &  \cdots & 0 \\
    0   &  0  & \ddots &  0 \\
    0   & \cdots & 0 &  K(\T\, E_{N-1})-O(\T\, E_{N-1})
  \end{matrix}
\right),
\end{equation*}
and therefore
\begin{eqnarray*}
\|K(\T\, \H)
- O(\T\, \H) \| \leq
\|\mathcal{E}\| =
 \max_{j=0,1,\ldots, N-1} \|K(\T\, E_j) -O(\T\, E_j)\|.
\end{eqnarray*}
Since $|E_j| \leq \beta$, $j=0,1,\ldots,N-1$, we finally arrive at
\begin{equation}
\label{eq:err_splitting}
  \|(q_1 + i \, p_1) -\e^{-i\,\T\, \H} (q_0+i\, p_0)\| \leq
\epsilon(\beta \, \T) \,
\|q_0 + i\, p_0 \|,
\end{equation}
where
\begin{equation}
\label{eq:epsilon}
 \epsilon(\theta) =
 \sup_{-\theta\leq y \leq \theta} \|K(y) -O(y)\|.
\end{equation}

 By taking into account that $\det K(y) \equiv 1$, the 2-norm of the $2\times2$ matrix $K(y)-O(y)$ can be explicitly computed to give
 \begin{eqnarray*}
\|K(y) -O(y)\| &=&
     \sqrt{(C(y)-\cos(y))^2+(S(y)-\sin(y))^2} \\
& & + \sqrt{C(y)^2+S(y)^2-1},
 \end{eqnarray*}
where
\begin{equation}
  \label{eq:CS}
  C(y) = \frac12 (K_{11}(y) + K_{22}(y)), \quad
S(y) = \frac12 (K_{12}(y) - K_{21}(y)).
\end{equation}
Notice that $\det K(y) \equiv 1$ implies
\begin{equation*}
  C(y)^2+S(y)^2-1 = \frac14 (K_{11}(y) - K_{22}(y))^2 + \frac14 (K_{12}(y) + K_{21}(y))^2
\end{equation*}
and thus $C(y)^2+S(y)^2-1\geq 0$ for all real values of $y$.

\subsubsection{Error estimates for several steps of a splitting method}
\label{sss:err_n}

Ideally, given a positive integer $m$ and $\theta=\beta\, \T>0$, one would like to determine a sequence (\ref{eq:sequence}) of real numbers so that
$\epsilon(\theta)$ is minimized. 
The error bound $\epsilon(\theta)$ being small implies that the $(2m)$th degree polynomial $C(y)$ 
(resp. the $(2m+1)$th degree polynomial $S(y)$) is
a good polynomial approximation of $\cos(y)$ (resp. $\sin(y)$) for $y\in [-\theta,\theta]$, which implies that increasingly large
values of $\theta=\beta\, \T$ will require longer sequences of coefficients (that is, larger values of $m$), and consequently more computational work. The
situation here is in complete analogy with what happened to Taylor and Chebyshev polynomial approximations in the previous section.

By applying the methodology exposed in \cite{blanes11eao} we have
determined several sequences (\ref{eq:sequence}) of length $2m+1$ of (near-to-optimal) coefficients for $m$ up to $60$. The procedure 
is described in detail in the Appendix. As shown there, the task is by no means trivial, and severe technical difficulties arise 
when trying to extend the procedure to
arbitrarily large values of
$\theta=\beta\, \T$ (and hence arbitrarily long sequences of coefficients). This is in contrast with Taylor and Chebyshev approximations.

This drawback can always be circumvented  by approximating the solution $z(\T) = O(\T\, \H)z_0$ of the system of
ordinary differential equations (\ref{td.1b}) in
the standard step-by-step way. In our case, approximating $z(\T)$ in $n$ steps of length
\begin{equation*}
  \tau=\frac{\T}{n}
\end{equation*}
simply consists in approximating $O(\T\, \H)z_0 = O(n\, \tau\, \H)z_0 = O(\tau\, \H)^nz_0$ by the vector $K(\tau\, \H)^nz_0$,
where $K(y)$ is a $2\times2$ matrix with polynomial entries (defined in terms of the sequence (\ref{eq:sequence}) as before)
that should approximate the rotation matrix $O(y)$ for $y \in [-\frac{\beta\, \T}{n},\frac{\beta\, \T}{n}]$.

Clearly, the resulting procedure for approximating $\e^{-i\,\T\, \H}u_0$ can be written as an algorithm of the form (\ref{alg:splitting}), corresponding to a sequence of coefficients (with a $(2m)$-periodic pattern) of length $2nm+1$. The corresponding error can be estimated as
\begin{eqnarray}
\nonumber
  \|(q_1 + i \, p_1) -\e^{-i\,\T\, \H} (q_0+i\, p_0)\| &\leq&
\|K(\tau\H)^n - O(n \tau \H) \| \,
\|q_0 + i\, p_0\|\\
  \label{eq:errS}
&\leq &
\epsilon^{(n)}(\tau\beta) \,
\|q_0 + i\, p_0 \|,
\end{eqnarray}
where
\begin{equation*}
  \epsilon^{(n)}(\theta) =
 \sup_{-\theta\leq y \leq \theta}\|K(y)^n -O(n\, y)\|.
\end{equation*}
Our goal is then to minimize $\epsilon^{(n)}(\theta)$. A reasonable requirement is that  $K(y)^n$ be bounded for all $n$.
This only happens
in general for a certain range of values of $y$. One thus defines the stability threshold $y_*$ as the largest non
negative real number such that $K(y)^n$ is bounded independently of $n \ge 1$ for all $y \in (-y_*,y_*)$ \cite{blanes08otl}.
In particular, for the sequence (\ref{eq:abStrang}) corresponding to the application of $m$ steps of the
Strang splitting,  the stability threshold is $ y_* = 2m$. As a matter of fact, $2m$ is precisely the maximal stability threshold a sequence of coefficients (\ref{eq:sequence}) of length $2m+1$ can achieve \cite{jeltsch81soe}.

From the analysis carried out in \cite{blanes11eao}, it is possible to show that
 \begin{eqnarray*}
   \|K(y)^n -O(n\, y)\| &\leq& 2\sin(n (\arccos(C(y))-y)/2) \\
&& +
\sqrt{\frac{S(y)^2}{1-C(y)^2}-1} + \frac{1}{2}
\left(
\frac{S(y)^2}{1-C(y)^2}-1
\right),
 \end{eqnarray*}
provided that  $y \in [-y_*,y_*]$.
This implies that, if $\tau\beta \leq y_*$, then
\begin{eqnarray}
\nonumber
  \|K(\tau\H)^n -O(n\tau\H)\|  &\leq &
 \sup_{-\tau\beta\leq y \leq \tau\beta} \|K(y)^n -O(n\, y)\|  = \epsilon^{(n)}(\tau \beta) \\ 
  \label{eq:errS2}
&\leq& n\mu(\tau\beta)+\nu(\tau\beta),
\end{eqnarray}
where
\begin{eqnarray}
\label{eq:mu}
  \mu(\theta) & = & \sup_{-\theta\leq y \leq \theta} |\arccos(C(y))-y|,  \\
\label{eq:nu}
  \nu(\theta) & = &  \sup_{-\theta\leq y \leq \theta}
\sqrt{\frac{S(y)^2}{1-C(y)^2}-1} + \frac{1}{2}
\left(
\frac{S(y)^2}{1-C(y)^2}-1
\right).
\end{eqnarray}

As mentioned before, we have determined several optimized splitting methods of $m$ stages (determined by a sequence
of coefficients (\ref{eq:sequence}) of length $2m+1$) for $m$ up to $60$.
The relevant parameters of such splitting methods are collected in 
Table~\ref{tab:relevant_parameters}. In this table,
 $M_m^{(\gamma)}$ refers to a method of $m$ stages, with error coefficients $\epsilon(\theta)$, $\mu(\theta)$, $\nu(\theta)$ optimized for $\theta = \gamma m$. For instance, method 
 $M_{60}^{(1.3)}$ can be used to approximate $\e^{-i t H} u_0$ with an error bounded (according to (\ref{eq:err_splitting}) and Table~\ref{tab:relevant_parameters}) by $1.2 \times 10^{-9} \|u_0\|$ provided that $ |t| \leq 78/\beta$. Furthermore,  $\e^{-i t H} u_0$ can be approximated by applying $n$ steps of length $\tau=t/n \leq\tau_{\mathrm{max}}:=78/\beta$ of method 
$M_{60}^{(1.3)}$ with an error bounded (according to (\ref{eq:errS2})) by 
\begin{equation*}
(7.8 n \times 10^{-11} + 1.2 \times 10^{-9}) \|u_0\|.
\end{equation*}

 In some cases two methods with the same values of $m$ and $\gamma=\theta/m$ have been collected, in which case they are labeled $a$ and $b$.
For instance, methods $M_{60}^{(1.4)a}$ and $M_{60}^{(1.4)b}$ are both designed to 
approximate $\e^{-i t H} u_0$ with $n$ steps of length $\tau=t/n \leq \tau_{\mathrm{max}} :=84/\beta$ of the method. However, they differ in the actual error estimate (\ref{eq:errS2}): in the first case, the error is bounded (provided that $\beta |t| \leq 84 n$)  by 
$(2.4 n \times 10^{-8} + 7.4 \times 10^{-8}) \|u_0\|$, while the second one admits the error estimate $(3.7 n \times 10^{-9} + 2.6 \times 10^{-6}) \|u_0\|$. This means that method 
$M_{60}^{(1.4)a}$ will be more efficient if $\beta |t| \leq 10452$, and the opposite otherwise.

Thus, given the upper bound $\beta$ of the spectral radius of $\H$ and the target time $\T$, if one wants to approximate
  $\e^{-i\,\T\, \H}u_0$ by applying $n$ steps of method $M_{60}^{(1.4)a}$, one should choose the smallest positive integer $n$ such that
\[
  \frac{\T}{n} \leq \tau_{\mathrm{max}}:= \frac{84}{\beta},  \quad \mbox{ that is, } \quad  n=\mbox{\texttt{Ceiling}}[\T\beta/84].
\]
For instance, suppose the target time $\T$ and the bound $\beta$ are such that $\T\beta=1000$. Then, clearly, $n=12$, so that 12 steps of
scheme $M_{60}^{(1.4)a}$ have to be applied with step size $\tau = 1000/(12 \beta) \simeq 83.3/\beta$ to achieve the target time. In this way
one gets an approximation with estimated error of size $3.62\times 10^{-7}\, \|u_0\|$ with a computational work ($2nm=2\times 12\times 60=1440$ real matrix-vector products of the form $H v$) comparable to the use of a Chebyshev polynomial approximation of degree $720$. In contrast, to guarantee a similar precision with Chebyshev, a polynomial of degree at least $1135$ is required, since  this is the minimum value of $m$ such that $\epsilon_{m}^{C}(1000)\, \|u_0\| \leq 3.62\times 10^{-7}\, \|u_0\|$, with $\epsilon_{m}^{C}(\theta)$ given in (\ref{eq:ErrorChebyshev}).

It is worth remarking the error coefficients for Strang splitting method (\ref{eq:abStrang}) with the same value of
$\gamma=\theta/m = 1.4$ (also collected in Table \ref{tab:relevant_parameters}) are much larger than for methods $M_{60}^{(1.4)a}$ and $M_{60}^{(1.4)b}$.

\subsubsection{Error estimates for combined splitting methods}

Sometimes it is just more efficient to apply a combination of two different methods instead of $n$ steps of the same scheme. For instance,
suppose that $\T\beta=177$
and we have an error tolerance of \texttt{tol}=$10^{-7}$. If we use  $M_{60}^{(1.4)a}$ then
$\T\beta/84 \simeq 2.1$ so that, according with the previous considerations, method $M_{60}^{(1.4)a}$ has
to be used with $n=3$ steps of size $\tau = 59/\beta$, much smaller than the value $\tau_{\mathrm{max}} =84/\beta$
for which the scheme has been designed.
This would result in an approximation fulfilling the required error tolerance obtained with $360$ real matrix-vector products of the form $H v$. A better strategy would be
the following: apply two steps of scheme $M_{60}^{(1.4)a}$ with step size $\tau_{\mathrm{max}}= 84/\beta$ to approximate 
$w=\e^{-i\,  2\tau_{\mathrm{max}}\, \H}u_0$
 and then approximating $\e^{-i\,  (\T-2 \tau_{\mathrm{max}})\, \H}w$ by using some other method with less stages. More generally,
we take $n=\mbox{\texttt{Floor}}[\T\beta/84]$ steps of length $\tau_{\mathrm{max}} = 84/\beta$ to get $w=\e^{-i\,  n \tau_{\mathrm{max}} \, \H}u_0$
and then we approximate $\e^{-i\,  (\T- n \tau_{\mathrm{max}})\, \H}w$ with another method of Table \ref{tab:relevant_parameters} involving less stages.

To decide which method
has to be used for this last step, we need an error estimate for the approximation obtained with such a combination of two methods.
Assume that we apply $n$ steps of length $\hat \tau$ of a method characterized by a $2\times 2$ matrix $\hat K(y)$ with polynomial entries, followed by a step of length $\tau$ of a method characterized by the matrix $K(y)$, where $\T=n \hat \tau + \tau$. From the preceding considerations, it is enough to estimate  $\|K(\tau\, \H) \hat K(\hat \tau\, \H)^n - O(t\, \H)\|$. This can be done in terms of the functions $\hat \mu(\theta)$, $\hat \nu(\theta)$ associated to  $\hat K(y)$ as defined in subsection~\ref{sss:err_n}, and the error function $\epsilon(\theta)$ associated to $K(\theta)$ as in subsection~\ref{sss:err}, together with the following function associated to $K(y)$:
\begin{equation}
  \label{eq:delta}
  \delta(\theta) = \sup_{-\theta \leq y \leq \theta}\|K(y)\|-1.
\end{equation}
Indeed, one obtains the following error estimate:
\begin{eqnarray}
  \label{eq:error-comb-meth}
  \|K(\tau\, \H) \hat K(\hat \tau\, \H)^n - O(t\, \H)\| &\leq&
\|K(\tau\, \H) - \e^{-i\, \tau \H}\| \, \|O(n \hat \tau \H)\|    \nonumber \\
&& + \|K(\tau\, \H)\| \, \|\hat K(\hat \tau \, \H)^n - O(n \hat \tau \H)\| \\
&\leq&
\epsilon(\tau \beta) + (1+\delta(\tau \beta)) (n \, \hat \mu(\hat \tau\, \beta) + \hat \nu(\hat \tau \, \beta)).  \nonumber
\end{eqnarray}
Since, as it can be noticed in Table \ref{tab:relevant_parameters}, $\delta(\theta) \simeq \epsilon(\theta) \ll 1$, then we can take simply
\begin{equation}  \label{eq:error-comb-meth-simpl}
  \|K(\tau\, \H) \hat K(\hat \tau\, \H)^n - O(t\, \H)\| \lesssim  
\epsilon(\tau \, \beta) + n \, \hat \mu(\hat \tau\, \beta) + \hat \nu(\hat \tau \, \beta).
\end{equation}
It is worth remarking that such an approximation will require $2 (n \hat m + m)+1$ real matrix-vector products of the form $H v$, and thus is equivalent in complexity to the application of a  (Chebyshev or Taylor) polynomial approximation of degree $n \hat m + m$.

\begin{table}
\centering
\begin{eqnarray*}
  \begin{array}{|c||c|c|c|c|c|c|c|} \hline
 &   & \theta=  &   &  &  &  & \\
M_m^{(\theta/m)} &  m &
\beta \tau_{\mathrm{max}}  & y_*/m  & \epsilon(\theta) & \mu(\theta) & \nu(\theta) & \delta(\theta)\\
\hline\hline
M_{10}^{(0.5)} & 10 &   5   &0.63 & 3.6\times 10^{-8} & 8.7\times 10^{-11} &
   9.8\times 10^{-8} & 3.6\times 10^{-8}\\
 \hline
M_{10}^{(0.9)} & 10 &  9 & 0.94 & 3.4\times 10^{-5} & 2.9\times10^{-5} &1.1\times 10^{-5}& 6.0\times 10^{-6} \\
 \hline\hline
M_{20}^{(0.6)} & 20 &   12   & 0. 79 & 1.6\times 10^{-13} & 1.4\times 10^{-13} & 5.8\times 10^{-14} &   2.5\times 10^{-14}\\
 \hline
M_{20}^{(1)} & 20 &    20    &    1.1 & 4.1\times 10^{-7} &
   1.8\times 10^{-8} & 4.8\times 10^{-7} & 4.0\times 10^{-7}\\
 \hline\hline
M_{30}^{(0.75)} & 30 &  22.5   &  0.84 & 8.1\times 10^{-15} & 3.3\times 10^{-16} & 1.5\times 10^{-14} & 7.9\times 10^{-15} \\ \hline
M_{30}^{(1)} & 30 &    30    & 1.0 & 4.1\times 10^{-10} & 1.9\times
   10^{-10} & 3.1\times 10^{-10} &
   2.6\times 10^{-10} \\ \hline
M_{30}^{(1.3)} & 30 &  39 &  1.36  &  2.3\times 10^{-5} & 5.2\times
   10^{-6} & 2.2\times 10^{-5} &
   2.0\times 10^{-5} \\
 \hline\hline
  M_{40}^{(1)}& 40 &    40    &  1.1 & 1.8\times 10^{-12} & 4.9\times 10^{-14} & 2.4\times 10^{-12} & 1.8\times 10^{-12} \\
 \hline
M_{40}^{(1.2)} & 40 &    48    &  1.26 & 2.1\times 10^{-8} & 2.1\times 10^{-8} & 5.3\times 10^{-10} & 4.7\times 10^{-10} \\
 \hline
M_{40}^{(1.4)} & 40 &    56    &  1.48  &1.48 \times 10^{-5} & 4.0\times 10^{-6} & 1.7\times 10^{-5} & 1.7\times 10^{-5} \\
 \hline\hline
 M_{50}^{(1)}& 50 & 50 & 1.07 & 4.5\times 10^{-15} & 4.5\times 10^{-15} &
   2.0\times 10^{-17} & 1.8\times 10^{-17} \\
 \hline
    M_{50}^{(1.1)}& 50 & 55 & 1.13 & 4.5\times 10^{-13} & 4.2\times 10^{-13} &
   4.1\times 10^{-14} & 3.5\times 10^{-14} \\
 \hline
 M_{50}^{(1.2)}& 50 & 60 & 1.26 & 5.4\times 10^{-11} & 2.7\times 10^{-11} &
   3.8\times 10^{-11} & 3.4\times 10^{-11} \\ \hline
M_{50}^{(1.3)a} & 50 & 65 & 1.32 & 1.2\times 10^{-8} & 1.2\times 10^{-8} &
   8.3\times 10^{-10} & 7.6\times 10^{-10} \\
\hline
  M_{50}^{(1.3)b}& 50 & 65 & 1.32 & 5.9\times 10^{-7} & 9.5\times 10^{-11} &
   6.1\times 10^{-7} & 5.9\times 10^{-7} \\ \hline\hline
M_{60}^{(1.1)} & 60 &  66    & 1.15    & 7.2\times 10^{-15} & 7.2\times 10^{-15} &
   2.6\times 10^{-17} & 2.2\times 10^{-17}
 \\ \hline
M_{60}^{(1.2)a} & 60 &  72    & 1.3 & 1.5\times 10^{-12} & 1.1\times 10^{-12} & 8.3\times 10^{-13} & 7.5\times 10^{-13} \\ \hline
M_{60}^{(1.2)b} & 60 &  72    & 1.26 & 4.2\times 10^{-11} & 6.5\times 10^{-14} & 4.6\times 10^{-11} & 4.2\times 10^{-11} \\
 \hline
 M_{60}^{(1.3)} & 60 &  78    &   1.36 & 1.2\times 10^{-9} & 7.8\times 10^{-11} & 1.2\times 10^{-9} & 1.2\times 10^{-9} \\
 \hline
M_{60}^{(1.4)a} & 60 &  84    &     1.41 & 8.4\times 10^{-8} &
   2.4\times 10^{-8} & 7.4\times 10^{-8} & 7.1\times 10^{-8}\\
\hline 
M_{60}^{(1.4)b} & 60 &  84    &    1.46 & 2.9\times 10^{-6} & 3.7\times 10^{-9} & 2.9\times 10^{-6} & 2.9\times 10^{-6} \\
 \hline
\hline\hline
\mathrm{Strang} & 1 &  1    & 2 &  1.8\times 10^{-1} & 4.7\times 10^{-2} & 1.5\times 10^{-1} & 1.3\times 10^{-1} \\
\hline
\mathrm{Strang} & 1 &  1.4    & 2 &  5.1\times 10^{-1} & 1.5\times 10^{-1} & 4.0\times 10^{-1} & 4.0\times 10^{-1} \\
\hline
\mathrm{Strang} & 1 &  1.9   & 2 & 1.34862 & 0.606472 &  2.4894 & 1.1746   \\ \hline
  \end{array}
\end{eqnarray*}
 \caption{Relevant parameters of several symplectic splitting methods
especially designed to integrate the semi-discretized
Schr\"odinger equation using a  time step $\tau=t/n$ with a maximum value 
$\tau_{\mathrm{max}} $.  Here  $y_*$ stands for the stability threshold and
$\epsilon(\theta)$, $\mu(\theta)$, $\nu(\theta)$, and $\delta(\theta)$ (for $\theta= \beta \tau_{\mathrm{max}}$) are the coefficients (appearing in the
error estimates obtained in Subsection 3.2) given in (\ref{eq:epsilon}), (\ref{eq:errS2}), and (\ref{eq:delta}) respectively.
}
  \label{tab:relevant_parameters}
\end{table}

\subsection{Flow of the algorithm}
\label{flowalgo}

Once a set of symplectic splitting methods constructed for providing approximations under different
conditions are
available (methods collected in Table~\ref{tab:relevant_parameters}) we still have to design a strategy to
select the most appropriate scheme and step-size to carry out
the numerical integration in time with the desired accuracy and a as small as possible
computational cost.


The user has to provide the values for $\Emin$ and $\Emax$,
a subprogram to compute the product $Hv$ for a given real vector $v$, the final
integration time $\T$ and the desired error tolerance \texttt{tol}. The procedure then
implements the shifting (\ref{shifting}), computes the value of $\beta$ and determines the normalized Hamiltonian
$\H$.

Next, the algorithm determines the most efficient method (or composition of methods) among the list of available schemes 
which provides the desired result: it chooses the cheapest method with error bounds below such tolerance and, if several methods with the same computational cost (same value of $m$) satisfy this condition, the algorithm chooses the scheme with 
the smallest error bound. This can be achieved if one starts the search from the methods with 
the smallest value of $m$ and, for each value of $m$, proceeds by decreasing accuracy, i.e. by increasing the value of 
$\theta = \beta \tau_{\mathrm{max}}$. For a given value of $t\beta$ and \texttt{tol} the algorithm checks 
for each method if $t\beta\leq \beta \tau_{\mathrm{max}}$ and, if this condition is satisfied, then it examines if $\epsilon(\theta)<$\texttt{tol}. 
This procedure corresponds to the sequence of methods collected in Table~\ref{tab:relevant_parameters} from top to bottom.

If none of the methods from the table satisfy both conditions for $t\beta$ and \texttt{tol}, then the time integration is split, i.e. $t\beta$ is
divided and a composition of one or several methods is used instead. Due to the high performance of the methods with the largest number of stages (in this case 60) the algorithm examines the cost of $n$ steps for the six 60-stage methods where $n=\mbox{\texttt{Floor}}[t\beta/\tau_{\mathrm{max}}\beta]$ and the last step is carried using one method from the list of methods. It chooses the cheapest methods with the smaller error bound among the composition of methods which provide the desired accuracy.

In this way, if we denote by $K^{(\gamma)}_{m}$ the matrix associated to method $M^{(\gamma)}_{m}$, then  
the resulting splitting method corresponds to the composition
\begin{equation}\label{eq:composition}
    K^{(\gamma_2)}_{m}(\tau\beta) \left(\hat K^{(\gamma_1)}_{60}(\hat \tau\beta)\right)^{n_1},
\end{equation}
where the algorithm chooses the methods (labelled by
$\gamma_1,\gamma_2,m$), the time steps, $\tau,\hat \tau$, and the
value of $n_1$, where $n_1=0$ if the method uses just one step. If
$n_1>0$ the error bound is given by
(\ref{eq:error-comb-meth-simpl}) while for $n_1=0$ the error bound
is just $\epsilon(\tau \, \beta)$.

This strategy has been implemented as a
 Fortran code which is freely available for download at the website \cite{website15}, together with some notes and examples illustrating the
 whole procedure.

In order to compare the efficiency of the resulting algorithm with the polynomial approximations based in Taylor and Chebyshev, 
with the error estimates collected in Table~\ref{tab:relevant_parameters} 
we have represented in Figure \ref{fig:mTayChev} the computational work (equivalent to a polynomial approximation of degree $m$) required for different tolerances and values of $\beta t$.
Diamonds, squares and circles correspond to the error tolerances $10^{-4}$, $2 \cdot 10^{-7}$ and $10^{-11}$, respectively,
obtained with
one or several steps of schemes in Table~\ref{tab:relevant_parameters}.
Notice that our algorithm based on symplectic splitting  methods provide better accuracy
with a considerably reduced  computational effort.

\section{Numerical examples}
\label{sec.5}

Next we apply the algorithm based on symplectic splitting methods presented in section \ref{sec.3} to two different
examples and compare its main features with Chebyshev and Taylor polynomial approximations. For the first example,
previously considered in \cite{lubich08fqt} to illustrate Chebyshev and Lanczos
approximations, we provide in addition the codes we have produced to generate the results and figures collected here. These can be 
found at \cite{website15}.
The second example illustrates
the performance of the methods on a one-dimensional Schr\"odinger equation with a smooth potential.

\paragraph{Example 1.}

The problem consists in computing $u(t) = \exp(-i t \widetilde{H}) u_0$ with $u_0 \in \mathbb{C}^N$ a unitary random vector and
the tridiagonal matrix
\begin{equation}\label{eq:Ex1Matrix}
 \widetilde H = \frac{1}{2} \left( \begin{array}{ccccc}
   2 & -1 & & & \\
   -1 & 2 & -1 & &  \\
   & & \ddots & & \\
    & & -1 & 2 & -1 \\
   & & & -1 & 2
   \end{array}  \right) \in \mathbb{R}^{N\times N}.
\end{equation}
The eigenvalues of $\widetilde H$ verify $0\leq E_k \leq 2$ for all $k$, so that we can take
$\Emin=0$, $\Emax=2$, and thus $\alpha = \beta = 1$ in (\ref{shifting}). In consequence, the problem reduces to
approximate
\begin{equation}  \label{scaledham}
   \e^{-i\alpha t}   \e^{-i\beta t \H} u_0,
   \qquad \mbox{where} \qquad \H =\widetilde H-I.
\end{equation}
We take $N=10000$ for the numerical experiments, but the results
are largely independent of $N$ (this is so even for the simplest, scalar case $N=1$).


Both Chebyshev and Taylor methods have been implemented in such a way that only
real valued matrix-vector products are used (we always
separate into the real and imaginary parts, i.e. $\H
u=\H (q+i\, p)=\H q+i\, \H p)$), so that Chebyshev
requires to store only 7 real vectors instead of 4 complex
vectors.

\begin{figure}[tb]
 \centering
  \includegraphics[width = 1.1\textwidth]{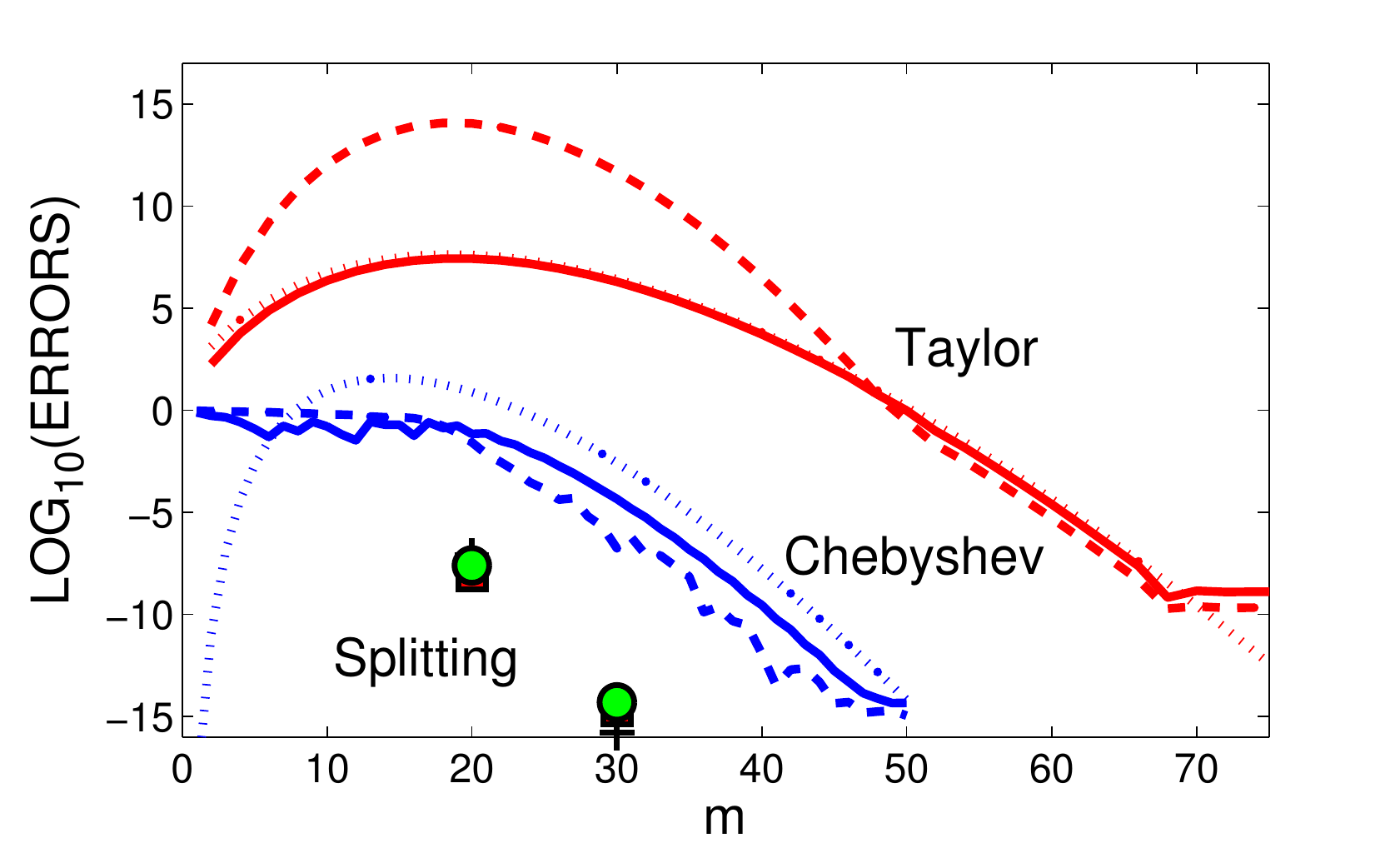}
\caption{{Different approximations to $\e^{-i \T \beta \H} u_0$, with $H$ given in (\ref{eq:Ex1Matrix})-(\ref{scaledham}), $u_0$ a random
vector, $\beta=1$ and
$\T \beta=20$ versus the degree of the polynomials, $m$. The
figure shows the relative error in energy (dashed lines), the
error in unitarity (solid lines) and error bounds (dotted lines)
for Chebyshev and Taylor methods. The results for the first two splitting methods with $\beta \tau_{\mathrm{max}}\geq \beta \T=20$, $M^{(1)}_{20}$ and  $M^{(0.75)}_{30}$,  are also shown:
relative error in energy
(filled squares), error in unitarity (filled squares) and error
bounds (crosses). }}
 \label{Fig0}
\end{figure}

We take as final time $\T = 20$ and measure the error in energy,
the error in the preservation of unitarity and the tolerance for
different values of $m$, the degree of the corresponding polynomials. The results are
shown in Figure~\ref{Fig0} with the following notation: dashed lines for the relative error in energy,
solid lines for the error in unitarity, and dotted lines for the theoretical error bounds of the approximate solutions.


From the figure it is clear that
the theoretical error bounds for the Taylor method are
quite accurate for this example (since the bounds for $E_{\mathrm{min}}$ and $E_{\mathrm{max}}$ are sharp) 
and that for the effective time-step $\tau \beta$ considered, the error is exceedingly large for $m$ below 
reaching the super linear convergence regime.
This is
not the case for the Chebyshev method (notice that the estimate (\ref{eq:ErrorChebyshev}) is valid
only for $m> \tau \beta$) since the coefficients
$c_k$ of the polynomial (\ref{alg:Clenshaw}) do not grow  as much as in Taylor.
%
We also depict the results achieved by
 the first two splitting methods with $\tau_{\mathrm{max}}\beta\geq \T\beta=20$,
$M^{(1)}_{20}$ and $M^{(0.75)}_{30}$. For these schemes the corresponding relative
error in energy is represented by filled squares, the error in unitarity by filled
circles and the error bounds by crosses.

\begin{figure}[tb]
 \centering
  \includegraphics[width = 1.1\textwidth]{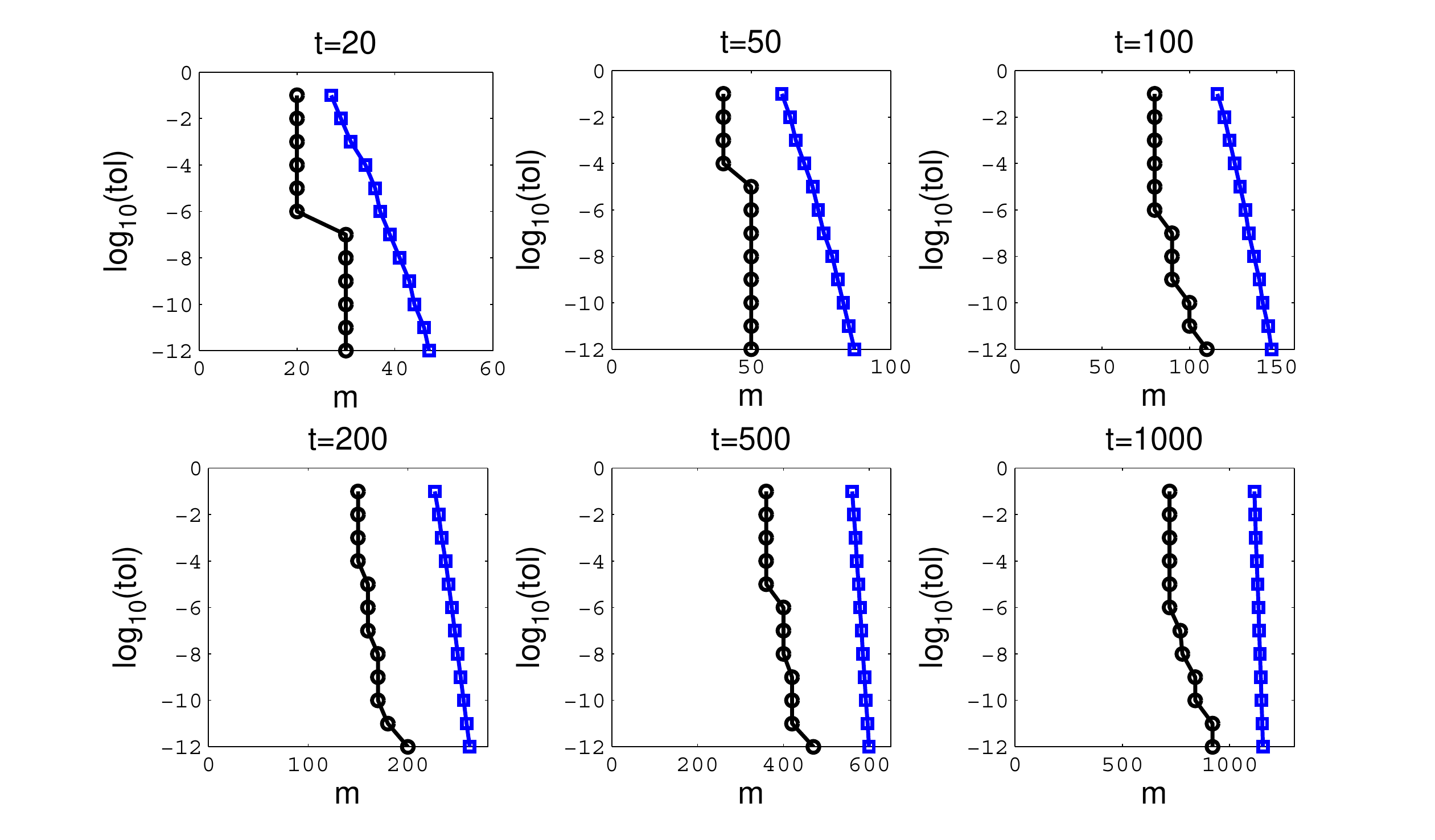}
\caption{{Degree $m$ of the polynomials to achieve tolerances
\texttt{tol}$\, =10^{-k}$, $k=1,2,\ldots,12$ for different values of $\beta \T $ ($\beta=1$ for this problem) as determined by the error
bound formulas using the Chebyshev method (squares) and the algorithm based on splitting methods (circles). }}
 \label{Fig1}
\end{figure}


The relative performance of different numerical integrators is usually tested by
measuring the error of the methods versus their computational
cost. However, the splitting methods we are considering in this work are
designed to achieve a given tolerance, whereas their
computational cost is determined through the error bound estimate. For
this reason, we believe it is more appropriate to measure the cost of
the methods for different values of the tolerance. In particular,
we take \texttt{tol}$\, =10^{-k}$, $k=1,2,\ldots,12$ and final integration times
$\T=20,50,100,200,500,1000$. Figure~\ref{Fig1} shows
the results obtained with Chebyshev (line with squares) and the algorithm based on splitting schemes (line with circles)
as a function of $m$. Even when
high accuracy is required over long integration times (the most advantageous situation for Chebyshev approximations),
the new algorithm requires a smaller value of $m$ and therefore
less computational effort. Notice how the algorithm selects the value of $m$ to achieve the desired tolerance.

\begin{figure}[tb]
 \centering
  \includegraphics[width = 1.1\textwidth]{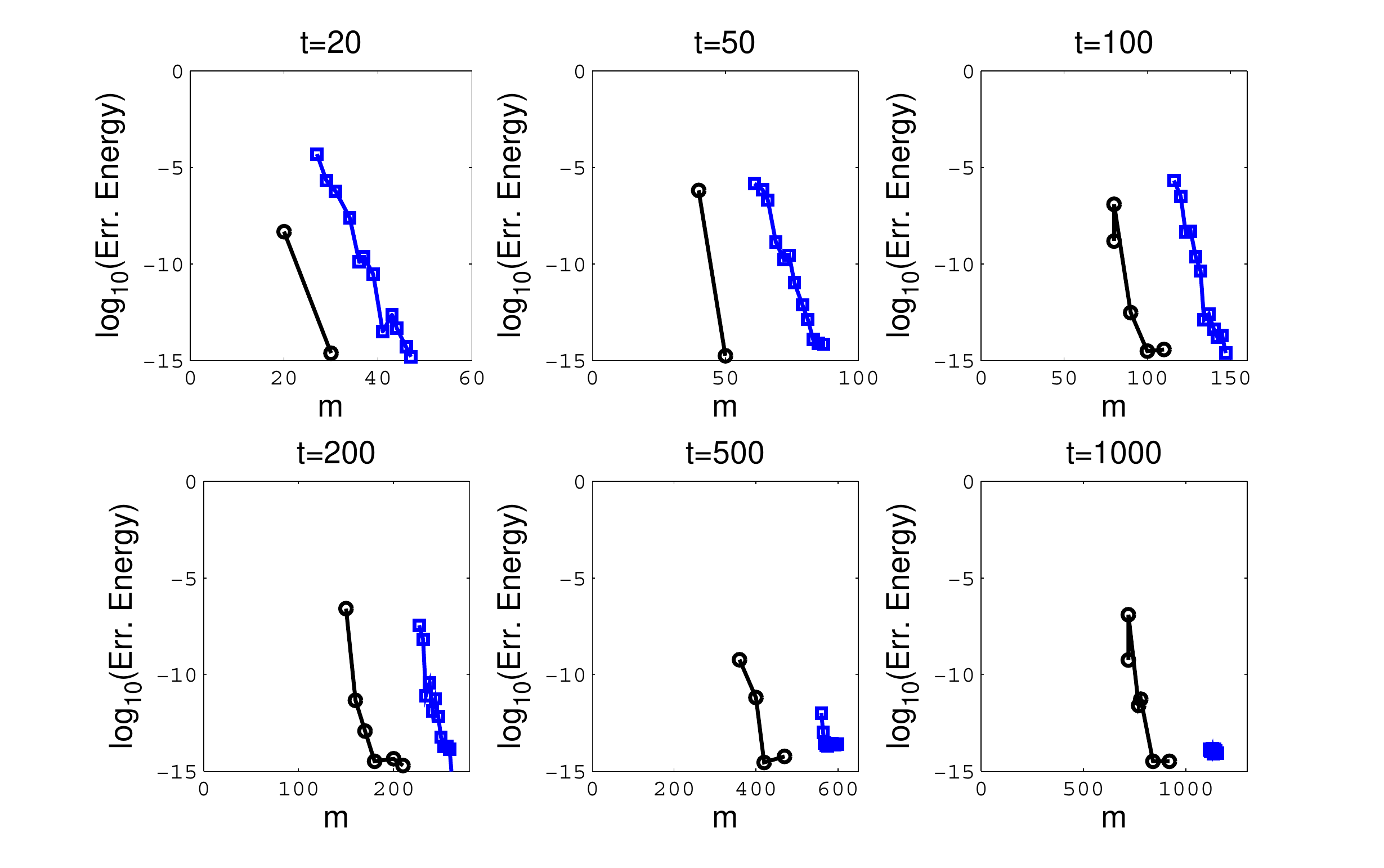}
\caption{{ Same as Figure~\ref{Fig1} but replacing the value of
the tolerance \texttt{tol} by the relative error in energy. }}
 \label{Fig1aa}
\end{figure}

Figure~\ref{Fig1aa} shows the corresponding results for the relative error in energy
versus $m$ for the same example.
Similar results
are obtained for the error in unitarity or the two-norm error for
which the error bounds apply (in this case one should compute
numerically the exact solution and compare with the
approximations obtained for each value of \texttt{tol}).

\paragraph{Example 2 (P\"oschl--Teller potential).}

To illustrate how the methods work on a more realistic case, we
consider the well known one-dimensional P\"oschl--Teller potential, which is an
anharmonic quantum potential
\[
  V(x) = -\frac{a^2}{2\mu}  \frac{\lambda(\lambda-1)}{\cosh^2(a x)},
 \]
with $a > 0, \lambda > 1$. It has
been frequently used in polyatomic molecular simulation and is
also of interest in supersymmetry, group symmetry, the study of
solitons, etc. \cite{shi07fmi,flugge71pqm,lemus02cot}. The
parameter $\lambda$ gives the depth of the well, whereas $a$
is related to the range of the potential. The energies are
 \[
 E_k =-\frac{a^2}{2\mu} (\lambda - 1 - k)^2, \qquad \mbox{ with }  \; 0\leq k \le \lambda - 1.
 \]

We take the following values for the parameters (in atomic units,
a.u.): reduced mass $\mu= 1745$ a.u., $a=2$,
$\lambda=24.5$ (leading to 24 bounded states), and $x\in[-5,5]$. Moreover, to apply a pseudo spectral space discretization we assume
periodicity of the potential in this range. The resulting $V(x)$ is thus
continuous and very close to differentiable for all $x \in \mathbb{R}$. Table~\ref{tab:2}
collects the bounds to the spectral radius (obtained according to
(\ref{eq:bounds})) and the corresponding shifting for the
P\"oschl--Teller potential when the space interval $x\in[-5,5]$
is split into $N$ parts and for different values of $N$. Notice how sensibly $E_{\mathrm{max}}$ depends on the space discretization.


\begin{table}[t]
\begin{eqnarray*}
  \begin{array}{|c|c|c|c|c|} \hline
    N & \Emin  & \Emax & \alpha & \beta  \\ \hline
   64 & -0.65988 & 0.11583 & -0.27202  & 0.38785 \\ \hline
  128 & -0.65988 & 0.46333 & -0.098275 & 0.5616  \\ \hline
  256 & -0.65988 & 1.8533  &  0.59672  & 1.2566  \\ \hline
  512 & -0.65988 & 7.4133  &  3.3767   & 4.0366 \\ \hline
 1024 & -0.65988 & 29.653  &  14.496   & 15.156  \\ \hline
  \end{array}
\end{eqnarray*}
 \caption{Bounds to the spectral radius and shifting for the  P\"oschl--Teller potential
 with the parameters considered in the text, when the space interval $x\in[-5,5]$ is split
 into $N$ parts.}
  \label{tab:2}
\end{table}

We take as initial condition a Gaussian function, $\psi
(x,0)=\sigma \, \e^{-(3x)^2}$, where $\sigma$ is a normalizing
constant, so the function and all its derivatives of practical
interest vanish up to round off accuracy at the boundaries. The
initial conditions contain part of the continuous spectrum, but
this fact is largely irrelevant due to the smoothness
of the periodic potential and wave function.


Suppose that one is interested in solving the corresponding semi discretized problem in time
with the following requirements:
\begin{itemize}
  \item[(I)] $N=128$, $\T = 15 \pi$, \texttt{tol}$\, = 10^{-9}$. In this case $\T \beta = 26.4648$.
  \item[(II)] $N= 512$, $\T = 40 \pi$, \texttt{tol}$\, = 10^{-6}$. Now $\T \beta = 507.254$.
\end{itemize}
 We have to determine first, of course, the degree $m$
of the polynomial  from the corresponding error bounds (for
Taylor the time interval is divided by two in (I) and by 36 in
(II) to avoid exceedingly large round off errors).
Table~\ref{tab:3} shows the number of matrix-vector products used
by each method (in bold) and the 2-norm error for each method
(compared with the exact solution obtained numerically with very
high accuracy).  In the first case our algorithm makes the
computations in a single step using $M_{30}^{(1)}$ while in the
second case it uses 6 steps of the scheme $M_{60}^{(1.4)a}$
followed by one step of $M_{10}^{(0.5)}$, i.e. the composition
(\ref{eq:composition}) is now
\[
    K^{(0.5)}_{10}(\tau\beta) \left(\hat K^{(1.4)a}_{60}(\hat \tau\beta)\right)^{6}
\]
with $\hat \tau=84/\beta$ and $\tau=40\pi-6\, \hat \tau$, and
for a total of 370 products. Again, the algorithm based on
symplectic splitting methods is able to produce results with the
required accuracy with less computational effort.


\begin{table}[t]
\begin{eqnarray*}
  \begin{array}{|c|c|c|c|}
\hline
   & Taylor & Chebyshev & Symplectic  \\ \hline
\hline  \begin{array}{c}
  {\T \, \beta=26.4648} \\ {\texttt{tol}=10^{-9}}
 \end{array}
 &
 \begin{array}{c}
  {\bf 104} \\
  3.4\times 10^{-12}
 \end{array}
   &
 \begin{array}{c}
  {\bf 51} \\
   3.7\times 10^{-12}
 \end{array}
    &
 \begin{array}{c}
   {\bf 30} \\
    4.2\times 10^{-11}
 \end{array}
      \\ \hline
 \begin{array}{c} {\T \, \beta=507.254} \\ {\texttt{tol}=10^{-6} }
 \end{array}   &
 \begin{array}{c}
  {\bf 1836} \\
  2.5\times 10^{-8}
 \end{array}      &
 \begin{array}{c}
  {\bf 587} \\
  3.4\times 10^{-15}
 \end{array}    &
 \begin{array}{c}
   {\bf 370} \\
    4.4\times 10^{-9}
 \end{array}
  \\ \hline
  \end{array}
\end{eqnarray*}
 \caption{Number of matrix-vector products (in bold) and actual errors
 given by the Taylor, Chebyshev and symplectic methods for different
$\T \, \beta$ and tolerances \texttt{tol}.}
  \label{tab:3}
\end{table}



\subsection*{Acknowledgements}

The authors acknowledge Ministerio de Econom\'{\i}a y Competitividad (Spain) for financial support through the coordinated project
MTM2013-46553-C3. AM is additionally partially supported by the Basque Government  (Consolidated Research Group IT649-13), and FC by NPRP GRANT \#5-674-1-114 from the Qatar National Research Fund.

\bibliographystyle{plain}

\

\section*{Appendix: Construction of methods}

We next describe the algorithm used to determine the coefficients (\ref{eq:sequence}) of length $2m+1$ for given $m$ and $\theta \in (0,2m)$.

Since all the error estimates in Subsection~\ref{sec.3.2} depend exclusively on the even polynomial (of degree $2m$)  
$C(y)$ and the odd polynomial (of degree $2m+1$) $S(y)$ given in (\ref{eq:CS}),  we first try to determine an appropriate pair of such polynomials 
satisfying the necessary conditions $C(0)=1$ and $C(x)^2+S(x)^2-1>0$ (for all $x\in \R$). Such pair of polynomials is uniquely determined 
by a polynomial $P(y)=C(y)+S(y)$ of degree $2m+1$ satisfying 
\begin{equation}
\label{eq:Pcond}
 P(0)=1, \qquad \frac{1}{2} (P(y)^2+P(-y)^2)-1\geq 0.
\end{equation}

Once an appropriate polynomial $P(y)=C(y)+S(y)$ satisfying (\ref{eq:Pcond}) is chosen, there is only a finite number of  corresponding sequences (\ref{eq:sequence}), which can be effectively determined~\cite{blanes08otl}. Since all of them share the same error estimates, we choose among them a sequence that minimizes
\begin{equation*}
\sum_{j=1}^{m+1} |a_j| +
\sum_{j=1}^{m} |b_j|.
\end{equation*}

We next focus on the effective construction of the polynomial  $P(y)=C(y)+S(y)$ of degree $2m+1$.

On the one hand, in order that the expression $\sqrt{C(y)^2+S(y)^2-1}$ featuring in the error estimate (\ref{eq:epsilon}) be small in the interval $y \in [-\theta,\theta]$,
\begin{equation}
\label{eq:csP}
\sup_{-\theta \leq y \leq \theta}
|\cos(y +e (y)) + \sin(y + e(y)) - P(y)|
\end{equation}
should be small for some real valued function $e(y)$. On the other hand,  minimizing
\begin{equation*}
\sqrt{(C(y)-\cos(y))^2+(S(y)-\sin(y))^2)}
\end{equation*}
in the interval $y \in [-\theta,\theta]$ is, provided that (\ref{eq:csP}) is small enough, essentially equivalent to minimizing
\begin{equation}
\sup_{-\theta \leq y \leq \theta} |e(y)|.
\end{equation}
To reduce the complexity of the final algorithm for determining the polynomial $P(y)$, we will try to minimize instead an 
alternative norm of $e(y)$ that we introduce next. First observe that if
\begin{equation}
\label{eq:e_cheb}
e(y) = \hat e_0 + \sum_{j\geq 1} \hat e_j \, T_{j}(y/\theta)
\end{equation}
 is the Chebyshev series expansion of the function $e(y)$, then
 \begin{equation}
 \label{eq:e_cheb_norm}
\sup_{-\theta \leq y \leq \theta} |e(y)| \leq  \sum_{j\geq 0} |\hat e_j |.
\end{equation}
This suggests that the right hand side of (\ref{eq:e_cheb_norm}) may be a good alternative to the supremum norm for sufficiently smooth functions $e(y)$. For practical considerations, we will minimize instead the following alternative norm of the function $e(y)$
\begin{equation}
\label{eq:norm}
\|e\|_{\theta} \equiv \sqrt{\sum_{j\geq 0} (\hat e_j) ^2}.
\end{equation}

 Now, to determine  the polynomial $P(y)=C(y)+S(y)$ of degree $2m+1$, we consider, for a given odd integer $l$ such that $m+1\leq l \leq 2m$, a given set of nodes $y_1,\ldots,y_{l}$ symmetrically placed in the interval $[-\theta,\theta]$, and a given odd polynomial $e(y)$ of degree $l-2$,  the polynomial $P(y)$ of degree $2l-1$ interpolating in the Hermite sense the function $\cos(y +e (y)) + \sin(y + e(y))$ for the nodes $y_1,\ldots,y_{l}$. In particular, this implies that $P(0)=1$ and
\begin{equation}
\label{eq:VW}
C(y)^2+S(y)^2-1 = \frac{1}{2} (P(y)^2+P(-y)^2)-1 = V(y) W(y)^2
\end{equation}
where $W(y)=(y-y_1)\cdots (y-y_l)$, and $V(y)$ is an even polynomial of degree $4m-2l+2$.
Thus,  $P(y)$ satisfies the necessary condition (\ref{eq:Pcond}) if and only if $V(y)\geq 0$ for all $y$.

  Notice that  the interpolation error (\ref{eq:csP}) admits an upper bound of the form
 \begin{equation}
\label{eq:csPW}
\sup_{-\theta \leq y \leq \theta}
|\cos(y +e (y)) + \sin(y + e(y)) - P(y)| \leq \frac{\eta}{(2l)!} \sup_{-\theta \leq y \leq \theta} W(y)^2,
\end{equation}
where  $\eta>0$ is an upper bound of the (absolute value of) the (2l)th derivative of the function  $\cos(y +e (y)) + \sin(y + e(y))$ in the interval $y \in [-\theta,\theta]$.

   For a prescribed set of nodes
$y_1,\ldots,y_l$,  we restrict the choice of the odd polynomial $e(y)$ (of degree $l-2$) so that the Hermite interpolating polynomial $P(y)$ is of degree $2m+1$ (which introduces $2(l-m)-2$ non-linear constraints on the non-zero coefficients $\hat e_1,\hat e_3, \ldots,\hat e_l$ of the polynomial $e(y)$ given by (\ref{eq:e_cheb})), and
determine $e(y)$ by minimizing the norm $\|e\|_{\theta}$ for that restricted set of odd polynomials $e(y)$ of degree $l-2$.
This produces a polynomial $P(y)$ for each choice of  the set of nodes $y_1,\ldots,y_l$.  It then remains to choose, for a prescribed positive odd integer $l$, such a set of nodes $y_1,\ldots,y_l$.

  The error estimate (\ref{eq:csPW}) suggests that a good choice for the interpolating nodes $\{y_1,\ldots,y_l\}$ may be given by the zeros of the Chebyshev polynomial $T_l(y/\theta)$ of degree $l$, which corresponds to minimizing the supremum norm (in the interval $[-\theta,\theta]$) of the polynomial $W(y)$. Notice that minimizing  the alternative norm $\|W\|_{\theta}$ also gives rise to the same set of nodes. It then only remains, for given odd positive number $2m+1$ and for given $\theta >0$, to determine the number $l$ of interpolating nodes, that should satisfy $m+1 \leq l \leq 2m$. If $l$ is too close to $2m$, then, very few degrees of freedom are left to minimize $\|e\|_{\theta}$, and if $l$ is too close to $m+1$, then the Hermite interpolating error (\ref{eq:csP}) is too large, causing the norm of the function $C(y)^2+S(y)^2-1$ not being small enough, in addition to $V(y)$ in (\ref{eq:VW}) typically not being positive. We thus proceed by determining $P(y)=C(y)+S(y)$ for different values of $l$ close to $(3m+3)/2$, and choosing, among those satisfying $V(y) \geq 0$,   one having the best error coefficient $\epsilon(\theta)$ defined in (\ref{eq:epsilon}).

  Unfortunately, choosing the interpolating nodes $\{y_1,\ldots,y_l\}$ as the zeros of the Chebyshev polynomial $T_l(y/\theta)$ of degree $l$ typically results in a polynomial $P(y)=C(y)+S(y)$ that does not satisfy the stability condition
  \begin{equation}
  \label{eq:stab_cond}
    |C(y)|\leq 1, \quad y \in [-\theta,\theta],
  \end{equation}
so that the error coefficients $\mu(\theta), \nu(\theta)$ are not well defined,
and thus the resulting splitting method cannot be reliably used in a step-by-step manner for large values of $t\beta$.
In order to produce splitting methods satisfying that stability condition for given $\theta$, we proceed iteratively to choose the
 interpolating nodes $\{y_1,\ldots,y_l\}$ and the corresponding polynomial $P(y)$ as follows: As a first approximation, we require the set of nodes $\{y_1,\ldots,y_l\}$ to contain the set $\{j \pi\ : \ j \in \mathbb{Z}, \ |j \pi|\leq \theta\}$ and determine the remaining nodes by minimizing the norm $\|W\|_{\theta}$ of $W(y)=(y-y_1) \cdots (y-y_l)$. Once the polynomial $P(y)=C(y)+S(y)$ is determined for that set
 of nodes $\{y_1,\ldots,y_l\}$, we compute the set of zeros of $C'(y))=0$ that are included in the interval $[-\theta,\theta]$ (that are typically close to $\{j \pi\ : \ j \in \mathbb{Z}, \ |j \pi|\leq \theta\}$), and determine the remaining nodes by minimizing the norm $\|W\|_{\theta}$ of $W(y)=(y-y_1) \cdots (y-y_l)$. 
 Successive iteration of this process gives a sequence of polynomials $P(y)=C(y)+S(y)$ that converge to a polynomial satisfying the stability condition (\ref{eq:stab_cond}).

\begin{figure}[htb]
 \centering
  \includegraphics[width = 0.9\textwidth]{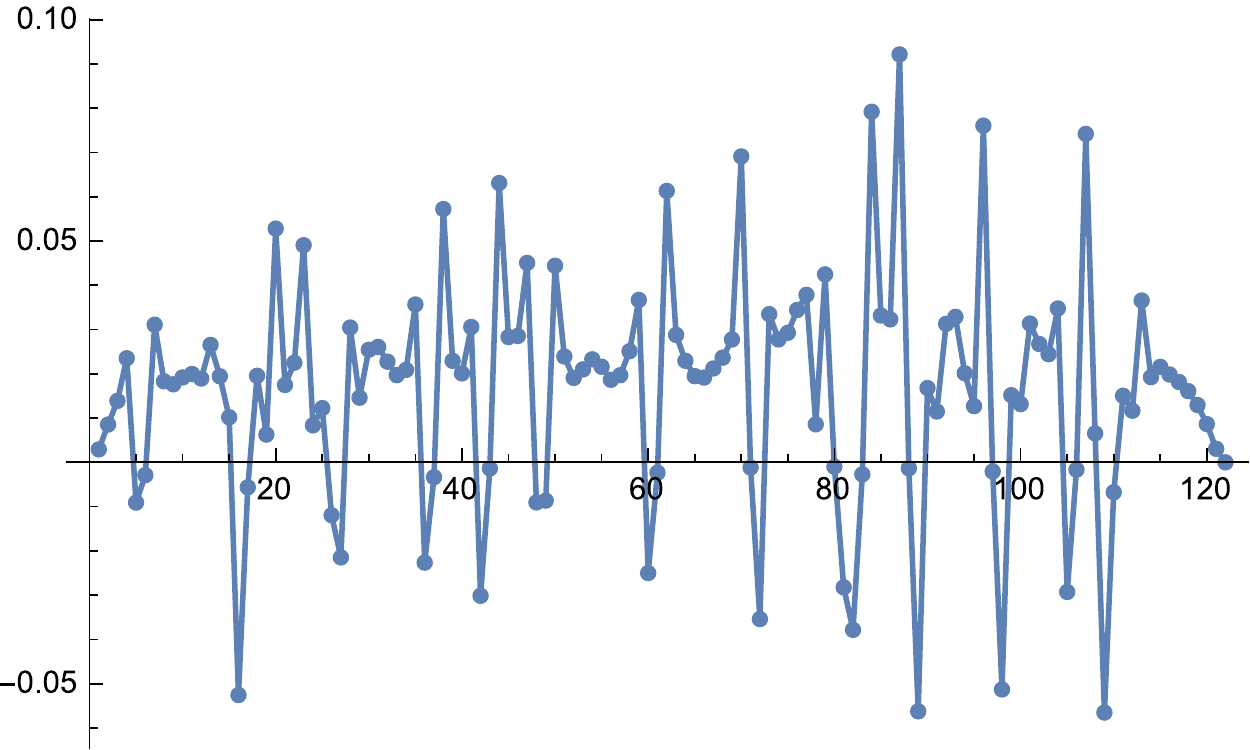}
\caption{Graphical representation of sequence $(a_1,b_1,a_2,b_2,\ldots,a_{60},b_{60},a_{61})$ of method $M_{60}^{(1.4)a}$ in Table~\ref{tab:relevant_parameters}, obtained with $\theta=84$ and $l=97$.
 \label{fig:coefs60Theta14}}
\end{figure}

 As an example, we have obtained the method $M_{60}^{(1.4)a}$ in Table~\ref{tab:relevant_parameters} by following this procedure for $m=60$, $\theta=84$, and $l=97$,  which has produced a splitting methods with sequence of coefficients (\ref{eq:sequence}) plotted in Figure~\ref{fig:coefs60Theta14}.

\end{document}